\documentclass[12pt]{article}
\usepackage{amsmath,amssymb,amsthm,amscd}
\usepackage{algorithm,algorithmic}
\usepackage[colorlinks, linkcolor=blue, anchorcolor=gree, citecolor=red]{hyperref}
\usepackage[numbers,sort&compress]{natbib}
\usepackage{graphicx}

\begin{document}

\newcommand{\Cyc}{{\rm{Cyc}}}\newcommand{\diam}{{\rm{diam}}}

\newtheorem{thm}{Theorem}[section]
\newtheorem{pro}[thm]{Proposition}
\newtheorem{lem}[thm]{Lemma}
\newtheorem{fac}[thm]{Fact}
\newtheorem{ob}[thm]{Observation}
\newtheorem{cor}[thm]{Corollary}
\theoremstyle{definition}
\newtheorem{ex}[thm]{Example}
\newtheorem{df}[thm]{Definition}
\newtheorem{remark}[thm]{Remark}

\renewcommand{\thefootnote}{}

\def\beql#1{\begin{equation}\label{#1}}
\title{\large\bf Lambda number of the power graph of a finite group}

\author{{Xuanlong Ma$^{1}$,~~Min Feng$^{2,}$\footnote{Corresponding author.}~,~~Kaishun Wang$^3$
}\\[15pt]
{\small\em $^1$School of Science, Xi'an Shiyou University, Xi'an, 710065, China}\\
{\small\em $^2$School of Science, Nanjing University of Science and Technology, Nanjing, 210094, China}\\
{\small\em $^3$Sch. Math. Sci. {\rm \&} Lab. Math. Com. Sys., Beijing Normal University, Beijing, 100875, China}\\
}

 \date{}

\maketitle

\begin{abstract}
The power graph $\Gamma_G$ of a finite group $G$ is the graph with the vertex set $G$, where two distinct elements are adjacent if one is a power of the other.
An $L(2, 1)$-labeling of a graph $\Gamma$ is an assignment of labels from  nonnegative integers to all vertices of $\Gamma$ such that vertices at distance two get different labels and adjacent vertices get labels that are at least $2$ apart.
The lambda number of $\Gamma$, denoted by $\lambda(\Gamma)$, is the minimum span over all $L(2, 1)$-labelings of $\Gamma$.
In this paper,
we  obtain bounds for $\lambda(\Gamma_G)$, and give necessary and sufficient conditions when the bounds are attained.
As applications, we compute the exact value of
$\lambda(\Gamma_G)$ if $G$ is a dihedral group,
a generalized quaternion group,  a $\mathcal{P}$-group or a cyclic group of order $pq^n$, where
$p$ and $q$ are distinct primes and $n$ is a positive integer.
\end{abstract}



\footnote{E-mail addresses: xuanlma@xsyu.edu.cn (X. Ma), fgmn\_1998@163.com (M. Feng), \\
wangks@bnu.edu.cn (K. Wang).}
\section{Introduction}
All graphs considered in this paper are finite, simple and undirected.
Let $\Gamma$ be a graph with the vertex set $V(\Gamma)$ and the edge set $E(\Gamma)$. The {\em distance} between vertices $x$ and $y$ is the length of a shortest path from $x$ to $y$ in $\Gamma$.
For nonnegative integers $j$ and $k$, an $L(j, k)$-labeling of $\Gamma$ is a nonnegative integer valued function $f$ on $V(\Gamma)$ such that $|f(u)-f(v)|\ge k$ whenever $u$ and $v$ are vertices of distance two and
$|f(u)-f(v)|\ge j$ whenever $u$ and $v$ are adjacent.
The {\em span} of $f$ is the difference between the maximum  and  minimum values of $f$. The {\em $L( j, k)$-labeling number} $\lambda_{j,k}(\Gamma)$ of $\Gamma$ is the minimum span over
all $L(j, k)$-labelings of $\Gamma$. The classical work of the $L(j, k)$-labeling problem is when $j = 2$ and $k = 1$. The $L( 2, 1)$-labeling number
of a graph $\Gamma$ is also called the {\em $\lambda$-number} of $\Gamma$ and denoted by $\lambda(\Gamma)$.

The problem of studying $L(j, k)$-labelings of a graph is motivated by the radio channel assignment problem \cite{Ha} and by the study of the scalability of optical networks \cite{Rob}.
In 1992, Griggs and Yeh \cite{GY92} formally introduced the notion of the
$L(j, k)$-labeling of a graph, and showed that the $L(2,1)$-labeling problem is NP-complete for general graphs.
 The $L(j, k)$-labelling problem, in particular in the $L(2, 1)$ case, has been studied extensively; see \cite{Ge,Ge2,W1,KiM} for examples. Surveys of results and open questions on the $L( j, k)$-labeling problem can be found in \cite{Ye}.

Graphs associated with groups and other algebraic structures have been actively investigated, since they have valuable applications
(cf. \cite{KeR}) and are related to automata theory (cf. \cite{K,K1}).
Zhou \cite{Zh1} studied
 $L(j,k)$-labelings of Cayley graphs of  abelian groups. Kelarev, Ras and Zhou \cite{Kerz} established connections between the structure of a semigroup and the minimum spans of distance labellings of its Cayley graphs. In this paper we study $L( 2, 1)$-labelings  of the {\em power graph} of a finite group.

The {\em undirected power graph} $\Gamma_G$ of a finite group $G$ has the vertex set $G$ and two distinct elements are adjacent if one is a power of the other. The concepts of a power graph
and an undirected power graph were first introduced by Kelarev and Quinn \cite{n1} and by Chakrabarty et al.  \cite{CGS}, respectively.
Since this paper deals only with undirected graphs, we use the term ``power graph" to refer to an undirected power graph.
Many interesting results on power graphs have been obtained in \cite{Cam,CGh,FMW,FMW1,kel2,MF,MFW,man}. A detailed list of results and open questions on power graphs can be found in \cite{AKC}.

Section~2 gives some preliminary results.
In Section~3, we obtain a sharp lower bound for the $\lambda$-number
of the power graph of a finite group $G$; as applications, we compute
$\lambda(\Gamma_G)$ if $G$ is a dihedral group,
a generalized quaternion group or  a $\mathcal{P}$-group.
In Section~4, we construct an upper bound for $\lambda(\Gamma_G)$, and classify all groups such that the upper bound is attained.

\section{Preliminaries}

A {\em path covering} of a graph $\Gamma$, denoted by $C(\Gamma)$, is a collection of vertex-disjoint paths in $\Gamma$ such that each vertex in $V(\Gamma)$ is contained in a path in $C(\Gamma)$. The {\em path covering number} $c(\Gamma)$ of $\Gamma$ is the minimum cardinality of a path covering of $\Gamma$. Let $\Gamma^c$  denote  the
complement of $\Gamma$.

\begin{pro}{\rm (\cite[Theorem 1.1]{Ge})}\label{Ge1}
Let $\Gamma$ be a graph of order $n$.

{\rm (i)} Then $\lambda(\Gamma)\le n-1$ if and only if $c(\Gamma^c)=1$.

{\rm (ii)} Let $r$ be an integer at least $2$.
Then $\lambda(\Gamma)=n+r-2$ if and only if $c(\Gamma^c)=r$.
\end{pro}

A vertex $x$ is a {\em cut vertex} in a grph $\Gamma$ if $\Gamma-x$ contains more connected components than $\Gamma$ does, where $\Gamma-x$ is the graph obtained by deleting the vertex $x$ from $\Gamma$.

\begin{pro}\label{prop1}
Let $\Gamma$ be a graph of order $n$ with a cut vertex $x$. Suppose that
all  connected components of $\Gamma-x$ are
$\Gamma_1,\ldots,\Gamma_t$ and $|V(\Gamma_i)|=n_i$ for $i\in\{1,\ldots,t\}$, where $n_t\le n_{t-1}\le \cdots \le n_1$. If
$
n_1\le \sum_{i=2}^{t}n_i,
$
then $\lambda(\Gamma)\le n$.
\end{pro}
\proof Note that $n\ge 3$. If  $n=3$, then $\Gamma$ is a path, and so $\lambda(\Gamma)=3$. In the following, suppose $n\ge 4$.
Write $\Delta=(\Gamma-x)^c$, and pick $v\in V(\Delta)$. Assume $v\in V(\Gamma_k)$ for some $k\in\{1,\ldots,t\}$. Then
$$
\deg_\Delta(v)\ge n-1-n_k,
$$
where $\deg_\Delta(v)$ is the degree of $v$ in $\Delta$. Since $\sum_{i=1}^tn_i=n-1$, we have $n_k\le n_1\le\frac{n-1}{2}$, which implies that $\deg_\Delta(v)\ge\frac{n-1}{2}$. It follows from Dirac's theorem (\cite[Theorem 4.3]{Bondy}) that $\Delta$ has a Hamilton cycle, and so $c(\Delta)=1$. Note that $\Delta=\Gamma^c-x$.
Then $c(\Gamma^c)\le 2$. By Proposition~\ref{Ge1}, we get the desired result.
\qed

\medskip

Let $\Gamma$ be a graph. A subset of $V(\Gamma)$ is
a {\em clique} if any two distinct vertices in this subset are adjacent
in $\Gamma$. The {\em clique number} $\omega(\Gamma)$ is the maximum cardinality of a clique in $\Gamma$. It is easy to see that
  \begin{equation}\label{s4}
  \lambda(\Gamma)\geq 2\omega(\Gamma)-2.
  \end{equation}
We give a sufficient condition for reaching  the lower bound in (\ref{s4}).

\begin{pro}\label{prop2}
  Let $C$ be a clique of a graph $\Gamma$ such that $|C|=\omega(\Gamma)$. Then $\lambda(\Gamma)=2\omega(\Gamma)-2$ if there exist partitions
  $$
  \{A_1,\ldots,A_s\}\qquad \text{and} \qquad \{C_1,\ldots,C_s,C_{s+1}\}
  $$
  of $V(\Gamma)\setminus C$ and $C$, respectively, satisfying the follows for each index $i\in\{1,\ldots,s\}$.

  {\rm(i)} $|A_i|\leq |C_i|-1$.

  {\rm(ii)} Every vertex in $A_i$ and every vertex in $C_i$ are nonadjacent in $\Gamma$.
\end{pro}
\proof
For  $1\le i\le s$ and $1\le j\le s+1$,
write
$$
A_i=\{u_{i1},\ldots,u_{im_i}\},\qquad C_j=\{v_{j1},\ldots,v_{jn_j}\}.
$$
Let $f$ be an integer valued function on $V(\Gamma)$ such that
$$
v_{1t}=2t,\quad f(v_{jt})=2(\sum_{k=1}^{j-1}n_k+t)\text{ for }j\geq 2, \quad f(u_{il})=f(v_{il})+1.
$$
By (i), we have $m_i\le n_i-1$, which implies that $f$ is well-defined,
and furthermore, the minimum and maximum values of $f$ are $2$ and $2\sum_{k=1}^{s+1}n_k$, respectively.
It follows from (ii) that $f$ is an $L(2,1)$-labeling of $\Gamma$, and so $$
\lambda(\Gamma)\le 2\sum_{k=1}^{s+1}n_k-2=2\omega(\Gamma)-2.
$$
By (\ref{s4}), the desired result follows.
\qed

\medskip

For a graph $\Gamma$,  a subset of $V(\Gamma)$ is
an {\em independent set} if no two of which are adjacent, and
the {\em independence number} $\alpha(\Gamma)$ is the maximum cardinality of an independent set in $\Gamma$.

\begin{pro}\label{prop3}
  Let $\Gamma$ be a graph of order $n$. Then
  $\lambda(\Gamma)\leq 2n-\alpha(\Gamma)-1$.
\end{pro}
\proof Write $r=c(\Gamma^c)$. If $r=1$, by Proposition~\ref{Ge1}~(i),
one has
$$
\lambda(\Gamma)\leq n-1\le 2n-\alpha(\Gamma)-1.
$$
Suppose $r\ge 2$. Let $A$ be an independent set of $\Gamma$ with $|A|=\alpha(\Gamma)$. Then  the subgraph induced by $A$ of $\Gamma^c$ has a Hamilton path, and so $c(\Gamma^c)\le n-|A|+1$. From Proposition~\ref{Ge1}~(ii), the desired inequality holds.
\qed

\medskip

The independence number of the power graph of a finite abelian group
has been studied. Denote by $\mathbb{Z}_n$ the cyclic group of order $n$.

\begin{lem}{\rm (\cite[Theorem 10]{Che})}\label{abin2}
Let $G$ be a finite abelian group. Then $\alpha(\Gamma_G)=2$ if and only if $G\cong \mathbb{Z}_{pq^n}$, where $p$ and $q$ are distinct
primes and $n$ is a positive integer.
\end{lem}

\begin{lem}{\rm (\cite[Theorem 5.4.10 (ii)]{Gor})}\label{pgroup0}
A $p$-group having a unique subgroup of order $p$ is either cyclic or
generalized quaternion.
\end{lem}

We extend Lemma~\ref{abin2} to the following result.

\begin{pro}\label{allgin2}
Let $G$ be a finite group. Then $\alpha(\Gamma_G)=2$ if and only if
$G\cong \mathbb{Z}_{pq^n}$, where $p$ and $q$ are distinct
primes and $n$ is a positive integer.
\end{pro}
\proof
The sufficiency follows from  Lemma~\ref{abin2}. Now suppose that $\alpha(\Gamma_G)=2$. Then $G$ has at most two
distinct subgroups of prime order, and furthermore the order of $G$ has at most two distinct prime divisors.

Suppose that $G$ is a $p$-group.
If $G$ has a unique subgroup of order $p$, by Lemma~\ref{pgroup0} we conclude that $\Gamma_G$ is complete or has independence number at least $3$, a contradiction. It follows that $G$ has two distinct subgroups of order $p$. Note that the center  $Z$ of $G$ is nontrivial.
Choose a subgroup $A$ of order $p$ in $Z$, and another subgroup $B$ of order $p$ in $G$. Then $AB\cong \mathbb{Z}_p\times\mathbb{Z}_p$, which implies that $G$ has at least $3$ subgroups of order $p$,  a contradiction.

The above contradiction implies that the order of  $G$ has exactly two two distinct prime divisors, say, $p$ and $q$.
By Lemma~\ref{pgroup0}, all Sylow subgroups of $G$ are cyclic.
Let $P$ be a Sylow $p$-subgroup of $G$ and $Q$ a Sylow $q$-subgroup
of $G$.
If there exists a Sylow $p$-subgroup $P_1$ of $G$ such that $P\ne P_1$, then $\{a,b,c\}$ is an independent set of $\Gamma_G$, where
$\langle a\rangle=P$, $\langle b\rangle=P_1$ and $\langle c\rangle=Q$,
a contradiction.
As a result, $G$ has a unique Sylow $p$-subgroup, and so
$P$ is normal in $G$. Similarly, we conclude that $Q$ is also normal in $G$. Consequently, $G$ is abelian. By Lemma~\ref{abin2}, we get the desired result.
\qed

\section{Lower bound}\label{dengyu}

In this section, we give a lower bound for the $\lambda$-number of the power graph of a finite group, and compute the exact value of
$\lambda(\Gamma_G)$ if $G$ is a dihedral group,
a generalized quaternion group or  a $\mathcal{P}$-group.

\begin{thm}\label{low}
Let $G$ be a group of order $n$. Then $\lambda(\Gamma_G)\ge n$, with  equality if and only if $(\Gamma_G-e)^c$ contains a Hamilton path,
where $e$ is the identity of $G$.
\end{thm}
\proof
Suppose that $f$ is an $L(2,1)$-labeling of $\Gamma_G$
with labels in $\{0,1,2,\ldots,\lambda(\Gamma_G)\}$.
Since $\Gamma_G$ has diameter at most $2$, the labels of all vertices under $f$ are pairwise distinct, which implies that $f$ is injective, and so $\lambda(\Gamma_G)\ge n-1$.
If $\lambda(\Gamma_G)= n-1$, then $f$ is a bijection, and hence
there exists an element $x\in G\setminus\{e\}$ such that
$|f(x)-f(e)|=1$, which contradicts that $x$ and $e$ are adjacent in
$\Gamma_G$. Therefore, we have $\lambda(\Gamma_G)\ge n$.

Note that $e$ is adjacent to every other vertex in $\Gamma_G$. Then $e$ is an isolated vertex in $(\Gamma_G)^c$, which implies that $c((\Gamma_G)^c)=2$ if and only if $(\Gamma_G-e)^c$ contains a Hamilton path. Hence, we get the desired result from Proposition~\ref{Ge1}~(ii).
\qed

\medskip

It is hard
to get a further characterization of the groups $G$ satisfying $\lambda(\Gamma_G)=n$,
where $n$ is the order of $G$.
In the remaining of this section, we give some groups such that the $\lambda$-numbers of their power graphs reach the lower bound in Theorem~\ref{low}.

\begin{ex}\label{}
For $n\ge 3$, we have $\lambda(\Gamma_{D_{2n}})=2n$, where $D_{2n}$ is the dihedral group of order $2n$.
\end{ex}
\proof Suppose
$$
D_{2n}=\langle a,b: a^n=b^2=e, bab=a^{-1}\rangle.
$$
Write $B=\{b,ab,a^2b,\ldots,a^{n-1}b\}$. Then $D_{2n}=\langle a\rangle\cup B$.
Define an integer valued function  $f$ on $D_{2n}$ as
$$
\begin{array}{clllll}
&f(e)=0,&f(a)=3,&f(a^2)=5,&\ldots,&f(a^{n-1})=2(n-1)+1,\\
&f(b)=2,&f(ab)=4,&f(a^2b)=6,&\ldots,&f(a^{n-1}b)=2(n-1)+2.
\end{array}
$$
Note that each element of $B$ is an involution. Then
$$
E(\Gamma_{D_{2n}})=E(\Gamma_{\langle a\rangle})\cup\{\{e,x\}: x\in B\},
$$
which implies that $f$ is an $L(2,1)$-labeling of $\Gamma_{D_{2n}}$, and so $\lambda(\Gamma_{D_{2n}})\le 2n$.  Thus, by
Theorem~\ref{low} we get $\lambda(\Gamma_{D_{2n}})=2n$.
\qed

\begin{ex}\label{q4n}
For $n\ge 2$, we have
$$
\lambda(\Gamma_{Q_{4n}})=\left\{
                                  \begin{array}{ll}
                                    4n+1, & \hbox{if $n$ is a power of $2$,}\\
                                    4n, & \hbox{otherwise},
                                  \end{array}
                                \right.
$$
where $Q_{4n}$ is the generalized quaternion group of order $4n$.
\end{ex}
\proof
Suppose
\begin{equation*}
Q_{4n}=\langle x,y: x^n=y^2, x^{2n}=1, y^{-1}xy=x^{-1}\rangle.
\end{equation*}
Then $y^{-1}=x^ny$, $|x^iy|=4$
and $(x^iy)^{-1}=x^{2n-i}y$ for $i\in\{1,\ldots,n-1\}$. Therefore
\begin{eqnarray}
V(\Gamma_{Q_{4n}})&=& \{e,x,\ldots,x^{2n-1}\}\cup
(\bigcup_{i=0}^{n-1}\{x^iy,(x^iy)^{-1}\}),\notag\\
E(\Gamma_{Q_{4n}})&=&E(\Gamma_{\langle x\rangle})\cup\bigcup_{i=0}^{n-1}E(\Gamma_{\langle x^iy\rangle})\label{s6}.
\end{eqnarray}
The power graph $\Gamma_{Q_{4n}}$ is shown in Figure~\ref{gq4n}.

\begin{figure}[hptb]
  \centering
  \includegraphics[width=8cm]{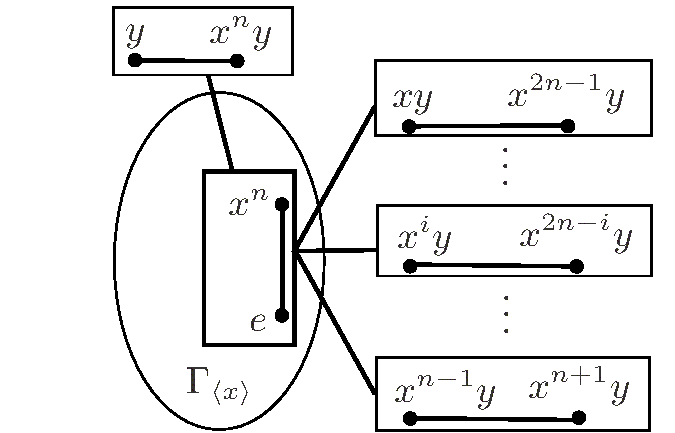}\\
  \caption{$\Gamma_{Q_{4n}}$}\label{gq4n}
\end{figure}

First suppose that $n$ is a power of $2$.
Then $x^n$ is adjacent to every other vertex of $\Gamma_{Q_{4n}}$, which implies that $x^n$ is an isolated vertex in $(\Gamma_{Q_{4n}}-e)^c$, and so $(\Gamma_{Q_{4n}}-e)^c$ does not have a Hamilton path.
According to Theorem~\ref{low}, one gets $\lambda(\Gamma_{Q_{4n}})\ge 4n+1$.
Define an integer valued function  $f$ on $V(\Gamma_{Q_{4n}})$ such that
$$
\begin{array}{l}
f(e)=0,f(x^n)=2,f(x)=2n,f(x^i)=2i \text{ for } i\in \{2,\ldots,2n-1\}\setminus\{n\},\\
f(x^jy)=2(j+1)+3 \text{ and } f(x^{2n-j}y)=2n+2(j+1)+1 \text{ for } j\in \{1,\ldots,n-1\},\\
f(y)=5,f(x^ny)=4n.
\end{array}
$$
It follows from (\ref{s6}) that $f$ is an $L(2,1)$-labeling of $\Gamma_{Q_{4n}}$ with span $4n+1$, and so $\lambda(\Gamma_{Q_{4n}})\leq 4n+1$. Thus, we have $\lambda(\Gamma_{Q_{4n}})=4n+1$.

In the following, suppose that $n$ is not a power of $2$. In view of Theorem~\ref{low}, it suffices to prove that $\lambda(\Gamma_{Q_{4n}})\le 4n$. Note that there exists an element $x_0$ of odd order in $\langle x\rangle$.
Write
$$\langle x\rangle\setminus\{e,x^n,x_0\}=\{z_1,z_2,\ldots,z_{2n-3}\}.$$
Let $f$ be an integer valued function on $V(\Gamma_{Q_{4n}})$ such that
$$
\begin{array}{l}
f(e)=0,f(x^n)=2,f(x_0)=3,f(z_i)=2(i+1)+1 \text{ for } i\in \{1,\ldots,2n-3\},\\
f(x^jy)=2(j+2) \text{ and } f(x^{2n-j}y)=2n+2(j+1) \text{ for } j\in \{1,\ldots,n-1\},\\
f(y)=4,f(x^ny)=4n-1.\\
\end{array}
$$
Since $x_0$ and $x^n$ are nonadjacent in $\Gamma_{Q_{4n}}$, it follows from (\ref{s6}) that $f$ is an $L(2,1)$-labeling of  $\Gamma_{Q_{4n}}$ with span $4n$, as desired.
\qed

\medskip

Now we give a sufficient condition for reaching  the lower bound in Theorem~\ref{low}.

\begin{lem}\label{co2}
Let $G$ be a noncyclic group of order $n$. Suppose that
all maximal cyclic subgroups of $G$ are $M_1,\ldots,M_t$, and
$|M_i|=n_i$ for $1\le i \le t$, where  $n_t\le n_{t-1}\le \cdots \le n_1$. If
\begin{equation}\label{s1}
|M_i\cap M_j|=1 \text{ for } 1\le i < j \le t
\end{equation}
and
\begin{equation}\label{s2}
n_1+t-2\le \sum_{i=2}^{t}n_i,
\end{equation}
then $\lambda(\Gamma_G)=n$.
\end{lem}
\proof Note that $e$ is a cut vertex of $\Gamma_G$.
 By (\ref{s1}), the connected components of $\Gamma_G-e$ are
 $\Gamma_{M_1}-e,\ldots,\Gamma_{M_t}-e.$
 Hence, the desired result  follows from Theorem~\ref{low} and Proposition \ref{prop1}. \qed

\medskip

A group is a {\em $\mathcal{P}$-group} \cite{CDLS} if every
nonidentity element of the group has prime order.

\begin{cor}\label{mainthm1}
Let $G$ be a $\mathcal{P}$-group of order $n$. Then
$$
\lambda(\Gamma_G)=\left\{
                                  \begin{array}{ll}
                                    2(n-1), & \hbox{if $n$ is a prime,}\\
                                    n, & \hbox{otherwise}.
                                  \end{array}
                                \right.
$$
\end{cor}
\proof
If $n$ is a prime, then $G$ is cyclic, and so $\Gamma_G$ is complete, which implies that $\lambda(\Gamma_G)=2(n-1)$. In the following, assume that $n$ is not a prime. Note that $G$ is noncyclic and each maximal cyclic subgroup of $G$ has prime order. With reference to Lemma~\ref{co2}, the equation (\ref{s1}) is valid. In order to get the desired result, it suffice to show that (\ref{s2}) holds.
Now by  \cite[Main Theorem]{CDLS}, one of the following cases occurs:

(a) $G$ is a $p$-group of exponent $p$, where $p$ is a prime.

(b) $G\cong A_5$, the alternating group on $5$ letters.

(c) $G$ is isomorphic to a Frobenius group $[P]Q$, where $P$ is a Sylow $p$-subgroup of $G$ and of exponent $p$, $Q$ is a Sylow $q$-subgroup of $G$ and of order $q$.

We use the notation in Lemma~\ref{co2} to verify (\ref{s2}). Note that if $n_1=n_2$ then (\ref{s2}) holds.
If (a) or (b) occurs, then  $n_1=n_2$. Suppose (c) occurs. Then $G$ has a unique Sylow $p$-subgroup $P$, and the number of Sylow $q$-subgroups is at least $p$. If $p<q$, then $n_1=n_2=q$.
If $p>q$ and $|P|>p$, then $n_1=n_2=p$. If $p>q$ and $|P|=p$, then $G$ has $p$ Sylow $q$-subgroups, which implies that
$$
t=p+1,\quad  n_1=p,\quad n_2=\cdots=n_t=q,
$$
 and so (\ref{s2}) holds.
\qed

\medskip

By Corollary~\ref{mainthm1}, we get the following example.

\begin{ex}\label{elementary abelian}
Let
$\mathbb{Z}_p^n$ denote the elementary abelian $p$-group of order $p^n$. Then $\lambda(\Gamma_{\mathbb{Z}_p^n})=p^n$ for $n\ge 2$.
\end{ex}

\section{Upper bound}

In this section, we shall prove the following result.

\begin{thm}\label{upper}
  Let $G$ be a group of order $n$.
  If $G$ is not cyclic of prime power order,
  then $\lambda(\Gamma_G)\leq 2n-4$, with equality if and only if
  $G$ is isomorphic to $\mathbb{Z}_2\times\mathbb{Z}_2$ or
$\mathbb{Z}_{2q}$, where $q$ is an odd prime.
\end{thm}

We remark that if $G$ is a cyclic group of prime power order then $\Gamma_G$ is complete and so its $\lambda$-number is equal to $2n-2$, where $n$ is the order of $G$. In order to prove Theorem~\ref{upper}, we first give some useful lemmas.

\begin{lem}\label{zpq}
Let $p$ and $q$ be distinct primes. For a positive integer $n$, we have
$$
\lambda(\Gamma_{\mathbb{Z}_{pq^n}})
=\left\{
\begin{array}{ll}
2q^{n-1}(pq-p+1)-2,&\text{if }p<q;\\
2q^{n}(p-1),&\text{if }q<p.
\end{array}\right.
$$
\end{lem}
\proof
For  $1\le i \le n$, write
\begin{eqnarray*}
&&X_i=\{g\in \mathbb{Z}_{pq^{n}}:|g|=pq^{n-i}\}
=\{x_{i1},\ldots,x_{ir_i}\},\\
&&Y_i=\{g\in \mathbb{Z}_{pq^n}:|g|=q^{n+1-i}\}
=\{y_{i1},\ldots,y_{is_i}\}.
\end{eqnarray*}
Let $Z=G\setminus(\bigcup_{i=1}^n(X_i\cup Y_i))$. Then
$$
Z=\{g\in \mathbb{Z}_{pq^n}:|g|=pq^{n}\text{ or }|g|=1\},
$$
$$
r_i=\varphi(pq^{i-1}),s_i=\varphi(q^i)\text{ and }|Z|=\varphi(pq^n)+1,
$$
where $\varphi$ is the Euler's totient function. Note that every vertex in $X_i$ and every vertex in $Y_i$ are nonadjacent in $\Gamma_{\mathbb Z_{pq^n}}$.

Suppose $p<q$. Then $(\bigcup_{i=1}^nY_i)\cup Z$ is a clique of order $\omega(\Gamma_{\mathbb Z_{pq^n}})$ in $\Gamma_{\mathbb Z_{pq^n}}$ by \cite[Theorem 2]{man}.
Note that $r_i\le s_i-1$.  By Proposition~\ref{prop2}, we have
$$
\lambda(\Gamma_{\mathbb{Z}_{pq^n}})
=2\omega(\Gamma_{\mathbb{Z}_{pq^n}})-2=2q^{n-1}(pq-p+1)-2.
$$

Suppose $q<p$. Then $(\bigcup_{i=1}^nX_i)\cup Z$ is a clique of order $\omega(\Gamma_{\mathbb Z_{pq^n}})$ in $\Gamma_{\mathbb Z_{pq^n}}$.
If $(p,q)\ne (3,2)$, then $r_i\le s_i-1$. It follows from Proposition~\ref{prop2} that
$$
\lambda(\Gamma_{\mathbb{Z}_{pq^n}})
=2\omega(\Gamma_{\mathbb{Z}_{pq^n}})-2=2q^{n}(p-1).
$$
In the following, assume $(p,q)=(3,2)$. Then $r_i=s_i=2^{n-i}$ for $1\le i\le n-1$ and $r_n=s_n+1=2$. Let $r_{n+1}=2^n+1$, and write
$$
Z=\{x_{n+1,1},\ldots,x_{n+1,r_{n+1}}\}.
$$
Define an integer valued function $f$ on $V(\Gamma_{{\mathbb{Z}}_{3\cdot2^n}})$ such that
$$
x_{1t}=2t,\quad f(x_{jt})=2(\sum_{k=1}^{n}r_k+t)\text{ for }j\geq 2, \quad f(y_{il})=f(x_{il})+1.
$$
For $1\le i\le n-1$, any vertex of $Y_i$ is not adjacent to any vertex of $X_i\cup X_{i+1}$ in $\Gamma_{{\mathbb{Z}}_{3\cdot2^n}}$.
Hence  $f$ is an $L(2,1)$-labeling of $\Gamma_{{\mathbb{Z}}_{3\cdot2^n}}$ with minimum value $2$ and maximum value $2\sum_{k=1}^{n+1}r_k$. Since
$$
2\omega(\Gamma_{{\mathbb{Z}}_{3\cdot2^n}})-2
=2^{n+2}=2\sum_{k=1}^{n+1}r_k-2,
$$
we get the desired result by (\ref{s4}).
\qed

\begin{lem}\label{lem}
  Let $G$ be a group of order $n$. If $G$ has pairwise distinct
  elements $u_1,u_2,u_3$ and $u_4$ such that $u_i$ and $u_{i+1}$
  are nonadjacent in $\Gamma_G$ for $1\le i\le 3$, then
   $\lambda(\Gamma)\leq 2n-5$.
\end{lem}
\proof Note that $n\ge 4$. If $c((\Gamma_G)^c)=1$,
then $\lambda(\Gamma)\leq n-1\leq 2n-5$
by Proposition~\ref{Ge1}~(i). Now suppose $c((\Gamma_G)^c)\geq 2$.
Since $(\Gamma_G)^c$ contains a path $(u_1,u_2,u_3,u_4)$, we have
$c((\Gamma_G)^c)\leq n-3$. It follows from
Proposition~\ref{Ge1}~(ii) that $\lambda(\Gamma)\leq 2n-5$.
\qed

\medskip

Now we prove Theorem~\ref{upper}.

\medskip

\noindent {\em Proof of Theorem~{\rm\ref{upper}}.} Since $\Gamma_G$ is not complete,
we have $\alpha(\Gamma_G)\geq 2$.
If $\alpha(\Gamma_G)=2$,
by Proposition~\ref{allgin2} and
Lemma~\ref{zpq}, we conclude that $\lambda(\Gamma_G)\leq 2n-4$,
with equality if and only if $G$ is isomorphic to
$\mathbb{Z}_{2q}$ for some odd prime $q$.
In the following, suppose $\alpha(\Gamma_G)\geq 3$.
From Proposition~\ref{prop3}, one has $\lambda(\Gamma_G)\leq 2n-4$.
If  $G$ is isomorphic to $\mathbb{Z}_2\times\mathbb{Z}_2$,
then $\lambda(\Gamma_G)=2n-4$ by Example~\ref{elementary abelian}.
Now assume $\lambda(\Gamma_G)=2n-4$. Then $\alpha(\Gamma_G)=3$
by Proposition~\ref{prop3}. It suffice to show that
$G$ is isomorphic to $\mathbb{Z}_2\times\mathbb{Z}_2$.

Suppose that $n$ is divisible by distinct primes $p$ and $q$,
where $p<q$. It follows from \cite[Section 4, I]{Fr} that
the number of subgroups of order $p$ and $q$ are $k_1p+1$
and $k_2q+1$, respectively, where $k_1$ and $k_2$ are
nonnegative integers.
If $k_1$ or $k_2$ is positive, then
$\Gamma_G$ has an independent set
containing $4$ pairwise distinct elements of prime order, contrary to
$\alpha(\Gamma_G)=3$. So $k_1=k_2=0$.
Let $P$ and $Q$ be a Sylow $p$-subgroup and a Sylow $q$-subgroup of $G$, respectively.
It follows from Lemma~\ref{pgroup0} that $P$ is either cyclic or generalized quaternion, and $Q$ is cyclic.
Pick an element $y$ of order $q$ in $Q$. If $p>2$, choose an element $z$ of order $p$ in $P$, then $y,z,y^{-1}$ and $z^{-1}$ are pairwise distinct elements in $G$ such that
\begin{equation}\label{s5}
  \{\{y,z\},\{z,y^{-1}\},\{y^{-1},z^{-1}\}\}\cap E(\Gamma_G)=\emptyset,
\end{equation}
and by Lemma~\ref{lem} we get $\lambda(\Gamma_G)\le 2n-5$, a contradiction. Hence $p=2$. If there exists an element $z$ of order $4$ in $P$, then (\ref{s5}) holds, a contradiction. Therefore,  $P$ is isomorphic to $\mathbb Z_2$, and so the generator $x$ of $P$ is the unique element of order $2$ in $G$, which implies that $x$ belongs to the center of $G$. If $Q$ has an element $z$ of order $q^2$, replace $y$ with $xy$ in (\ref{s5}), then we get a contradiction. Therefore $Q=\langle y\rangle$. If $n$ has a prime divisor $r$ with $r\notin\{p,q\}$, pick an element $z$ of order $r$ in $G$, then  (\ref{s5}) holds, a contradiction. It follows that $G$ is isomorphic to $\mathbb{Z}_{2q}$, and so Proposition~\ref{allgin2} implies that $\alpha(\Gamma_G)=2$, a contradiction.

The above contradiction implies that $G$ is a $p$-group. If $p\ge 3$, then by Lemma~\ref{pgroup0}, $G$ has at least $p+1$ subgroups of order $p$, contrary to $\alpha(\Gamma_G)=3$. Hence $p=2$.
If $G$ has a unique subgroup of order $2$, then Lemma~\ref{pgroup0}
and Example~\ref{q4n} give  $\lambda(\Gamma_G)=n+1$, a contradiction.
Thus, the number of subgroups of order $2$ in $G$ is at least $3$.
If $G$ has an element of order $4$, we obtain a contradiction by Lemma~\ref{lem}. It follows that $G$ is elementary abelian, and from
Example~\ref{elementary abelian} we deduce that $G$ is isomorphic to $\mathbb{Z}_2\times \mathbb{Z}_2$, as desired.
\qed

\section*{Acknowledgement}
This work was completed during Ma and Feng's visit to the Beijing Normal University.
Wang's research was supported by National Natural Science Foundation of China (11371204, 11671043) and the Fundamental Research Funds for the Central University of China.

\end{document}